\documentclass[10pt,letterpaper]{article}
\usepackage[top=0.85in,left=2.75in,footskip=0.75in]{geometry}

\usepackage{changepage}

\usepackage[utf8]{inputenc}

\usepackage{textcomp,marvosym}

\usepackage{fixltx2e}

\usepackage{amsmath,amssymb}

\usepackage{cite}

\usepackage{nameref,hyperref}

\usepackage[right]{lineno}

\usepackage{microtype}
\DisableLigatures[f]{encoding = *, family = * }

\usepackage{rotating}

\raggedright
\usepackage[aboveskip=1pt,labelfont=bf,labelsep=period,justification=raggedright,singlelinecheck=off]{caption}

\makeatletter
\renewcommand{\@biblabel}[1]{\quad#1.}
\makeatother

\date{}

\usepackage{lastpage,fancyhdr,graphicx}
\usepackage{epstopdf}
\fancyheadoffset[L]{2.25in}
\fancyfootoffset[L]{2.25in}
\renewcommand{\a}{\alpha}
\renewcommand{\b}{\beta}
\newcommand{\bea}{\begin{eqnarray}}
\newcommand{\eea}{\end{eqnarray}}
\newcommand{\f}[2]{\frac{#1}{#2}}
\newcommand{\eq}{&=&}
\newcommand{\nn}{\nonumber \\ }
\newcommand{\ve}{\varepsilon}
\renewcommand{\d}{\delta}
\newcommand{\area}{\int_{-\infty}^\infty }

\newcommand{\p}{\partial}
\newcommand{\pp}[2]{\f{\p #1}{\p #2}}
\newcommand{\s}{\sigma}

\newcommand{\siki}[1]{Eq. (\ref{#1})}
\newcommand{\hyou}[1]{Fig. \ref{#1}}
\newcommand{\refauthor}[2]{#1 #2,}
\newcommand{\referencepaper}[7]
{#1 
#2. 
#3, 
#4, 
#5(#6), 
#7.
}

\newcommand{\referencebook}[4]
{#1 
#2. 
#3, 
#4. 
}

\begin{document}
\vspace*{0.35in}

\begin{flushleft}
{\Large
\textbf\newline{Property Safety Stock Policy for Correlated Commodities Based on Probability Inequality
}
}
\newline
\\
Takashi Shinzato\textsuperscript{1,\textpilcrow}
\\
\bigskip
\bf{1} Mori Arinori Center for Higher Education and Global Mobility,
Hitotsubashi University, Kunitachi, Tokyo, Japan.
\\
\bigskip

%
%

\textpilcrow The author contributed to this work.



* takashi.shinzato@r.hit-u.ac.jp
\end{flushleft}
\section*{Abstract}

Deriving the optimal safety stock quantity with which to meet customer satisfaction is one
of the most important topics in stock management. However, it is difficult to control the stock management of correlated
marketable merchandise when using an inventory control method that was developed under the assumption
that the demands are not correlated. For this, 
we propose a deterministic approach that uses a probability inequality to derive a reasonable safety stock for the case in which we know the correlation between various commodities. Moreover, over a given lead time, the relation
between the appropriate safety stock and the allowable stockout rate is analytically derived, and the potential of
our proposed procedure is validated by numerical experiments.



\section{Introduction}
Safety stock management is one of the most important
issues related to inventory control.  
A well-known method is commonly used to determine the appropriate inventory, but it is
based on the assumption that the demands for stock items are independent; thus,
it is not appropriate when the demands are correlated. 
Moreover, it cannot be applied when the statistical properties (mean and variance,) of the demand 
distribution are unknown.

In recent
years, 
various studies \cite{Fotopoulos,Takemoto,Takemoto2,Shinzato} have actively investigated
a novel approach to resolving this difficulty; it is based on a probability inequality, in particular,
the Markov inequality or Chebyshev inequality, 
and derives the relation between the quantity of
safety stock and the allowable stockout rate. 
For instance, Fotopoulos et al. discuss 
the relation between the safety stock quantity and the allowable
stockout rate; they use the Chebyshev inequality, since the actual
stockout rate becomes smaller than the average of the polynomial
of demand \cite{Fotopoulos}. Takemoto et al. \cite{Takemoto} and Takemoto and Arizono \cite{Takemoto2} have used
the Hoeffding inequality for the limited demand information, of which only the minimum and maximum values are known, and calculated the safety stock quantity which can
guarantee a given allowable stockout rate. 
In particular, Takemoto and Arizono \cite{Takemoto2} assessed the stockout rate for the supply of two types of electric power, which were represented by 
two uncorrelated random variables. 
Shinzato and Kaku \cite{Shinzato} derived
an expression for the relation between the safety stock quantity and the allowable 
stockout rate. They considered the case in which the time series correlation (trend) in the demand sequence of one of the commodities was known; this information was 
included by using the Chernoff inequality.
Although various studies have used a probability inequality to estimate the appropriate safety stock quantity, most analyses of 
safety stock have been for individual commodities or uncorrelated demands. However, in practice, stock management usually must address multiple commodities with correlated demands; this has not been sufficiently studied. 
In order to address this, we use the Chernoff inequality to derive the relation
between the safety stock quantity and the allowable stockout rate, and we perform experiments that validate the effectiveness of our proposed method.

This paper is organized as follows. In section 2, we introduce the Chernoff
inequality as a way to analytically estimate the stockout rate over the lead time. We then use this method to derive the relation between the safety stock quantity and the allowable stockout rate. In section 3, we use the statistical
properties of the known demand distributions of two correlated commodities and conduct numerical simulations to verify the validity of our proposed approach. In section 4, we summarize our conclusions and discuss areas of future work on this topic.

\section{Safety stock management and model setting}

In this section, we begin by discussing an existing safety stock management policy
that is based on the assumption that the distributions of the commodities are independent, and
formulate the associated problem. Then, for demand distributions that are not independent, we propose
a new safety stock management for a given safety stock quantity and an allowable stockout
rate; we use the Chernoff inequality to resolve the 
problems with the existing method \cite{Schmidt}.

\subsection{Safety stock management based on assumption
of independency}

This subsection discusses an existing safety stock management
system that is based on the assumption that the demands of single commodity are independently and identically distributed; it assumes that the statistical
signature of the demands is well known and that the market is stationary.
Here, the demands of the commodity are $D_t,(t=1, 2, \cdots, L)$, where $L$ is the lead time;
the distance from the average $E[D_t]=\mu$
is 
$X_t(=D_t-\mu),(t=1.2.\cdots,L)$, and we consider the case in which the demands are
independently and identically distributed (i.i.d.) with mean $0$
and variance $\s^2$.
Then, the safety stock quantity 
based on the previous method with respect to the allowable stockout rate $\d$, $SS_{\rm pre.}$,  is given as, 
\bea
SS_{\rm pre.}
\eq\sqrt{L\s^2}\times k,
\label{eq1}
\eea
where $k$ is the safety stock coefficient, and the allowable stockout rate, $\d$,  satisfies
\bea
\d\eq
\int_k^\infty
\f{dz}{\sqrt{2\pi}}e^{-\f{z^2}{2}}.\label{eq2}
\eea
Moreover, although the order point of this commodity is denoted by
$L\mu+SS_{\rm pre.}$, we will consider a safety stock quantity under uncertain demand. \siki{eq1} is mathematically guaranteed by the central limit theorem. Note that for a large lead time $L$,
the sum of the average demand gap
over the lead time, $SS=\sum_{t=1}^LX_t$, asymptotically follows
the normal distribution with mean 0 and
variance $L\s^2$, since the demands are uncorrelated.

However, since this method strongly requires that the demands $X_t$ be independent, 
it is difficult to manage the
safety stock of commodities for which the time series of the demand is correlated   \cite{Fotopoulos,Shinzato} or for which  
the demands of different commodities are correlated \cite{Silver}. Previous studies \cite{Kottas1,Kottas2,Ray1,Ray2,Bagchi1,Bagchi2,Van,Lau} have investigated qualitatively and quantitatively the nonlinear
effects in cases which it is not possible to assume independent
demands. 
Furthermore, Ray \cite{Ray2}, Lau and Wang \cite{Lau}, and Zhang et al. \cite{Zhang}
considered in detail the influence of the correlation which is
inherent in the demands. In particular, Shinzato and Kaku \cite{Shinzato} showed that when the time effect was relatively strong, the stock management
policy based on \siki{eq1} results in a stockout rate that is about $36$ times
the acceptable rate. 

In the next subsection, we consider a mathematical framework which can be used to evaluate the
safety stock quantity without requiring the assumption of independence of the demands; we then use this framework to propose a
novel method for safety stock management.
\subsection{Chernoff inequality\label{sec2.2}}

In this subsection, we discuss the use of the Chernoff inequality for deriving the safety
stock quantity without the assumption of that the demands are independent.
Consider a random variable
$D$ that takes a discrete value and has a density function $p(D)$, and where $D$ is not limited to a special distribution,
such as a Poisson or normal distribution)
Then, we have the following inequality  \cite{Schmidt}:
\bea
Pr[\eta\le D]
\le e^{-u\eta}E[e^{uD}],\label{eq3}
\eea
where $\eta$ is a constant,
$u>0$ is a control parameter, and 
$E[f(D)]$ is the expectation of 
$f(D)$ (here and in the following, $E[\cdots]$ denotes the expectation). 
This inequality is called the
Chernoff inequality, and it holds regardless of the
distribution of the random variable.
Without loss of generality, we assume the
average of the random variable $D$ to be $0$. 
In order to prove the Chernoff inequality, we begin by defining the following step function:
\bea
\Theta(Z)
\eq
\left\{
\begin{array}{ll}
1&0\le Z\\
0&{\rm otherwise}
\end{array}.
\right.
\label{eq4}
\eea
Using this step
function, the probability of $\eta\le D$ can be rewritten as
\bea
\label{eq5}
Pr[\eta\le D]
\eq
\int_\eta^\infty
dDp(D)\nn
\eq
\area dDp(D)\Theta(D-\eta)\nn
\eq E[\Theta(D-\eta)].
\eea
Then, 
\bea
\label{eq6}
\Theta(Z)&\le& e^{uZ}, \qquad(\forall u>0, \forall Z\in{\bf R}),
\eea
holds. Thus, the Chernoff inequality,
\bea
\label{eq7}
Pr[\eta\le D]\eq
E[\theta(D-\eta)]\nn
&\le& E[e^{u(D-\eta)}]\nn
\eq e^{-u\eta}E[e^{uD}],
\eea
is obtained. Since the Chernoff inequality in
\siki{eq3} holds for an arbitrary nonnegative $u>0$, there exists a
minimum of the right-hand side of \siki{eq3}:
\bea
\label{eq9}
Pr[\eta\le D]
&\le&
\mathop{\min}_{u>0}
\left[
e^{-u\eta}E[e^{uD}]
\right]\nn
\eq e^{-R(\eta)},
\eea
where the rate function $R(\eta)$ is defined as 
\bea
\label{eq10}
R(\eta)\eq\mathop{\max}_{u>0}
\left\{
u\eta-\phi(u)
\right\},\\
\label{eq11}
\phi(u)
\eq\log E[e^{uD}],
\eea
in which $\phi(u)$ is the cumulant generating function, which is defined as
the logarithm of the moment-generating function $E[e^{uD}]$ (see \ref{app-a1}).

Furthermore, we can use the cumulant generating function $\phi(u)$ to
assess the rate function $R(\eta)$ for any $\eta\in{\bf R}$. 
Note that the
cumulant generating function is convex with respect
to $u$, and from the definition in \siki{eq11}, we see that \siki{eq10} 
is the
Legendre transformation \cite{Shinzato,Cengel}, and thus $R(\eta)$ 
is also a
convex function of $\eta$. 
We now derive the relation between the safety stock
quantity and the allowable stockout rate by using the upper bound  $e^{-u\eta}E[e^{uD}]$, 
which is a tight upper bound on the probability 
$Pr[\eta\le D]$ of the stochastic event currently observed, 
that is, $\mathop{\min}_{u>0}e^{-u\eta}E[e^{uD}]=e^{-R(\eta)}$.

Note that we use the density function of a continuous random variable $D$ in the above proof of the Chernoff
inequality, although it also holds for a discrete
random variable. In addition, the probability
of $\eta\ge D$, that is, $Pr[\eta\ge D]$, is
\bea
\label{eq8}
Pr[\eta\ge D]&\le&
e^{-u\eta}
E[e^{uD}],
\eea
when $u<0$; the discrete case can be shown in a similar way.

As in the above discussion but with respect to the random variables $D_1,D_2$ 
and constants $\eta_1,\eta_2$,
 we consider the Chernoff inequality
for the probability of the inequalities $\eta_1\le D_1$ and 
$\eta_2\le D_2$, 
that is, $Pr[\eta_1\le D_1,\eta_2\le D_2]$, and we obtain
\bea
\label{eq12}
Pr[\eta_1\le D_1,\eta_2\le D_2]
\eq
E[\Theta(D_1-\eta_1)
\Theta(D_2-\eta_2)
]\nn
&\le&
e^{-u_1\eta_1-u_2\eta_2}
E[e^{u_1D_1+u_2D_2}],\nn
\eea
where $u_1,u_2>0$. 
In addition, we can write the Chernoff inequality
for the probability of inequalities for $N$ random variables
$\vec{D}=(D_1,\cdots,D_N)^{\rm T}\in{\bf R}^N$ 
and $N$ constants $\vec{\eta}=(\eta_1,\cdots,\eta_N)^{\rm T}\in{\bf 
R}^N$:
\bea
\label{eq13}
Pr[\vec{\eta}\le
\vec{D}
]
&\le&
e^{-\vec{u}^{\rm T}\vec{\eta}}
E[e^{\vec{u}^{\rm T}
\vec{D}
}],
\eea
where there are $N$ control parameters 
$\vec{u}=(u_1,\cdots,u_N)^{\rm T}\in{\bf R}^N,(u_i>0)$,
$\eta_i\le D_i$  holds at each
component, and ${\rm T}$ denotes the transpose of a vector or matrix. 
Note that we did not need to assume that the random variables are independent, and so this holds where the demands are correlated.

In similar way to the derivation of \siki{eq12}, we obtain the
Chernoff inequality for the
probability of inequalities $\eta\le D_1$ and $\eta_2\ge D_2$,
 that is,
$Pr[\eta_1\le D_1,\eta_2\ge D_2]$:
\bea
\label{eq14}
Pr[\eta_1\le D_1,\eta_2\ge D_2]
&\le&
e^{-u_1\eta_1-u_2\eta_2}
E[e^{u_1D_1+u_2D_2}],\nn
&&\qquad
(\forall u_1>0,\forall u_2<0);
\eea
for the probability of inequalities $\eta_1\ge D_1$ and $\eta_2\le D_2$,
\bea
\label{eq15}
Pr[\eta_1\ge D_1,\eta_2\le D_2]
&\le&
e^{-u_1\eta_1-u_2\eta_2}
E[e^{u_1D_1+u_2D_2}],\nn
&&\qquad
(\forall u_1<0,\forall u_2>0),
\eea
and for the probability of inequalities $\eta_1\ge D_1$ and $\eta_2\ge D_2$, 
\bea
\label{eq16}
Pr[\eta_1\ge D_1,\eta_2\ge D_2]
&\le&
e^{-u_1\eta_1-u_2\eta_2}
E[e^{u_1D_1+u_2D_2}],\nn
&&\qquad
(\forall u_1<0,\forall u_2<0).
\eea

Here, two points should be noted. First, although we consider the probability
satisfied with both of the conditions $\eta\le D_1$ and $\eta_2\le D_2$
in \siki{eq12}, in practice, it is only necessary to impose one of them:
\bea
\label{eq17}
&&Pr[(\eta_1\le D_1)
\cup
(\eta_2\le D_2)
]\nn
\eq
E[
\Theta(D_1-\eta_1)
+\Theta(D_2-\eta_2)
-\Theta(D_1-\eta_1)
\Theta(D_2-\eta_2)
]\nn
\eq
E[
\Theta(D_1-\eta_1)
+
\Theta(D_2-\eta_2)
\Theta(\eta_1-D_1)
].
\eea
When it is not possible to directly evaluate this probability, we can obtain a tighter upper bound on it
by using the Chernoff inequality in a similar way to that shown in \siki{eq12}. That is, the Chernoff
inequality can be evaluated as follows:
\bea
\label{eq18}
&&Pr[(\eta_1\le D_1)
\cup
(\eta_2\le D_2)
]\nn
&\le&
e^{-u_1\eta_1}
E[e^{u_1D_1}]
+
e^{-u_2\eta_2-u_3\eta_1}
E[e^{u_2D_2+u_3D_1}],\nn
&&
(\forall u_1>0,\forall u_2>0,\forall u_3<0).
\eea
In order to avoid duplication, in the discussion below, we will only consider $(\eta_1\le 
D_1)\cap(\eta_2\le D_2)=[\eta_1\le D_1,\eta_2\le D_2]$, were $D_i$ is demand of item $i$.

Second, the probability that the sum of two random variables, $D_1+D_2$, 
exceeds the constant $\eta$, $Pr[\eta\le D_1+D_2]=E[\Theta(D_1+D_2-\eta)]$,
satisfies the following Chernoff inequality:
\bea
\label{eq19-1}
Pr[\eta\le D_1+D_2]
&\le&
e^{-u\eta}
E[e^{uD_1+uD_2}],
\eea
where $u>0$. There are two ways to
interpret \siki{eq19-1}. 
First, as in the above argument, $D_1$ and $D_2$ represent the demands of the respective commodities during the same
period, and if substitution is allowed, then
$D_1+D_2$ 
is the total demand of the substitutable commodities and $Pr[\eta\le D_1+D_2]$ 
is the probability that the sum of the demands of
the substitutable commodities exceeds some
constant $\eta$  \cite{Takemoto2,Zhang}. See \ref{appendix:A} for a discussion of the management of the safety stock of substitutable commodities. The
second interpretation is that $D_i$ 
is the demand of a single item during period $i$. Then, 
$D_1+D_2$ represents the total demand during two periods, and $Pr[\eta\le D_1+D_2]$ 
is the probability that total demand exceeds 
$\eta$ \cite{Shinzato}. 
We will develop the latter interpretation
of \siki{eq19-1}, and we will determine, over a given lead time, 
the relation between the safety stock quantity and the allowable stockout rate; this will be discussed in detail in subsections \ref{sec2.4} and \ref{sec2.5}.

\subsection{Relationship between the Chernoff inequality and the stockout
rate\label{sec2.4}}
In this subsection, we will discuss the relationship between the sum
of demands that are generated over lead time $L$, that is, $\sum_{t=1}^LD_t$, and
the allowable stockout rate $\d$. Similar to what we did in \siki{eq19-1},
 we replace the
random variable $D$ and the constant $\eta$ 
in the Chernoff inequality
in \siki{eq3} with
$D=\sum_{t=1}^LD_t$ and 
$\eta=L\mu+L\ve$. Then, we obtain
\bea
\label{eq24}
Pr
\left[
L\mu+L\ve
\le
\sum_{t=1}^LD_t
\right]
&\le&
e^{-LR(\ve)},
\eea
where the rate function per degree $R(\ve)$ and the cumulant
generating function per degree $\phi(u)$ are defined by
\bea
\label{eq25}
R(\ve)
\eq\mathop{\max}_{u>0}
\left\{u(\mu+\ve)-\phi(u)\right\},\\
\label{eq26}
\phi(u)
\eq
\f{1}{L}
\log
E\left[
\exp\left(u
\sum_{t=1}^LD_t
\right)
\right],
\eea
in which $\mu=E[D_t]$ is the mean of the demand $D_t$. Note that on the right-hand side of \siki{eq24}, we have  $e^{-LR(\ve)}$, not $e^{-R(\ve)}$. To justify this, consider the following: if one tosses a fair coin $L$ times, then the number of
heads, $X$, follows the binomial distribution, since the
probability of $X=n\in{\bf Z}$, 
that is, $Pr[X=n]$, and/or the probability
of 
the principal part of the upper bound of the probability of the sum of $L$
random variables, is proportional to the $L$th power, that is, $e^{-LR(\ve)}$, not $e^{-R(\ve)}$.

Thus, $Pr[L\mu+L\ve\le\sum_{t=1}^LD_t]$
in \siki{eq24} describes the
probability that the sum of demands generated over lead time $L$, 
$\sum_{t=1}^LD_t$, is larger than the constant $\eta=L\mu+L\ve$, as the stockout rate. In the context of stock management, the argument in this probability, $L\ve$, can be regarded as the safety
stock quantity over lead time $L$.
Similarly, with respect to $N$ commodities, as in \siki{eq13},
\bea
\label{eq27}
Pr\left[
L\vec{\mu}+L\vec{\ve}
\le
\sum_{t=1}^L
\vec{D}_t
\right]
&\le&
e^{-LNR(\vec{\ve})},
\eea
where $\eta_i=L\mu_i+L\ve_i$, $D_{it}$ is the demand of
commodity $i$ at time $t$, $D_i=\sum_{t=1}^LD_{it}$, $E[D_{it}]=\mu_t$, and the rate function per degree $R(\vec{\ve})$ and cumulant generating
function per degree $\phi(\vec{u})$ are 
\bea
\label{eq28}
R(\vec{\ve})
\eq\mathop{\max}_{\vec{u}>0}
\left\{
\f{1}{N}
\sum_{i=1}^N
u_i(\mu_i+\ve_i)
-\phi(\vec{u})
\right\},\\
\label{eq29}
\phi(\vec{u})
\eq
\f{1}{NL}
\log E
\left[
\exp\left(
\sum_{i=1}^N
u_i
\sum_{t=1}^LD_{it}
\right)
\right].
\eea
Moreover, when the time sequences of the demands of commodity $i$ with lead
time $L$ are not correlated, one degree of the cumulant
generating function $\phi(\vec{u})$ in \siki{eq29} can be summarized as
$\phi(\vec{u})=
\f{1}{N}
\log E
\left[
\exp\left(\sum_{i=1}^Nu_iD_{it}\right)
\right]
$.
\subsection{
Relation between safety stock quantity and allowable
stockout rate
\label{sec2.5}
}
The right-hand side of the Chernoff inequality in \siki{eq24} is 
equal to the allowable stockout rate $\d$, where $\d=e^{-LR(\ve)}$, and 
the probability on the left-hand side can be regarded as the
practical stockout rate. The safety stock quantity $L\ve$ is thus 
\bea
\label{eq30-1}
L\ve\eq
L
R^{-1}
\left(-\f{1}{L}
\log \d
\right),
\eea
where 
$R^{-1}(\cdots)$ is 
the inverse function of rate function, $R(\ve)$. Note that 
this formula does not require the assumption of that the demands are i.i.d., 
and is described by the above-discussed probability theory. We obtain the 
relation between 
the safety stock quantity and the 
allowable stockout rate
for multiple 
correlated commodities, in a way similar to \siki{eq30-1}. We cannot derive accurate safety stocks from \siki{eq30-1} which are 
compatible with the allowable stockout rate $\d$ in principle. 
However, we can derive the allowable stockout rate $\d$ corresponding to the safety stocks of 
$N$ commodities, $L\vec{\ve}=(L\ve_1,\cdots,L\ve_N)\in{\bf R}^N$, 
from the formula $\d=e^{-LNR(\vec{\ve})}$.

As an application of the limits on the safety stock derived from the Chernoff inequality, we note the following: 
(1) When $\ve<0$, that is, when safety stock quantity $L\ve$ is negative, $R(\ve)=0$; thus, 
it does not provide a useful upper bound; 
(2) When the probability distribution of the demand is unknown, and we cannot assess the moment-generating function, 
we cannot directly use the Chernoff inequality for safety stock management. However, 
if we can analytically determine the moment-generating function, we can do so. For instance,  
with respect to the Weibull distribution $p(X,\a,\b)=\f{\a}{\b}\left(\f{X}{\b}\right)^{\a-1}
\exp\left[{-\left(\f{X}{\b}\right)^\a}\right]$ with known parameters $\a$ and $\b$, then 
since $E[e^{uX}]=\sum_{k=0}^\infty\f{u^k\b^k}{\Gamma(k+1)}\Gamma\left(\f{k}{\a}
+1\right)$, and for a 
log-normal distribution 
$p(X)=\f{1}{\sqrt{2\pi\s^2}X}\exp\left[{-\f{(\log 
X-\mu)^2}{2\s^2}}\right]$ for known parameters $\mu$ and $\s^2$, then $E[e^{uX}]=\sum_{k=0}^\infty\f{u^ke^{k\mu+\f{\s^2k^2}{2}}}{\Gamma(k+1)}$. 
Note that 
since a random variable which follows a Weibull distribution and/or log-normal distribution does not have a finite upper bound, 
the analytical approach based on the Hoeffding inequality for the safety stock, which was developed in \cite{Takemoto}, does not apply.

Moreover, under a stable market and with an unknown demand distribution, when $M$ is 
sufficiently large, if the demand sequence of dimension $L$, $(D_{1a},D_{2a},\cdots,D_{La})
(a=1,\cdots,M)$ is given, then
the estimation of the moment-generating function, 
$E\left[e^{u\sum_{t=1}^LD_t}\right]$,
can be replaced by 
$\f{1}{M}\sum_{a=1}^Me^{u\sum_{t=1}^LD_{ta}}$, 
in the Chebyshev inequality, for any $C(>0)$:
\bea
&&Pr\left[\left|E\left[e^{{u\sum_{t=1}^LD_t}}\right]-\f{1}{M}\sum_{a=1}^Me^{u\sum_{t=1}^LD_{ta}}\right|
\right.\nn
&&
\left.\ge C\sqrt{E\left[
\left(
e^{u\sum_{t=1}^LD_t}-
E\left[e^{u\sum_{t=1}^LD_t}\right]
\right)^2
\right]}
\right]\nn
&\le&\f{1}{MC^2}.
\eea
The 
moment-generating function (or its logarithm, the 
cumulant generating function) is thus estimated precisely (\hyou{fig1-est}).

\begin{figure}[tb]
\begin{center}
\includegraphics[width=1.0\hsize,height=0.6\hsize]{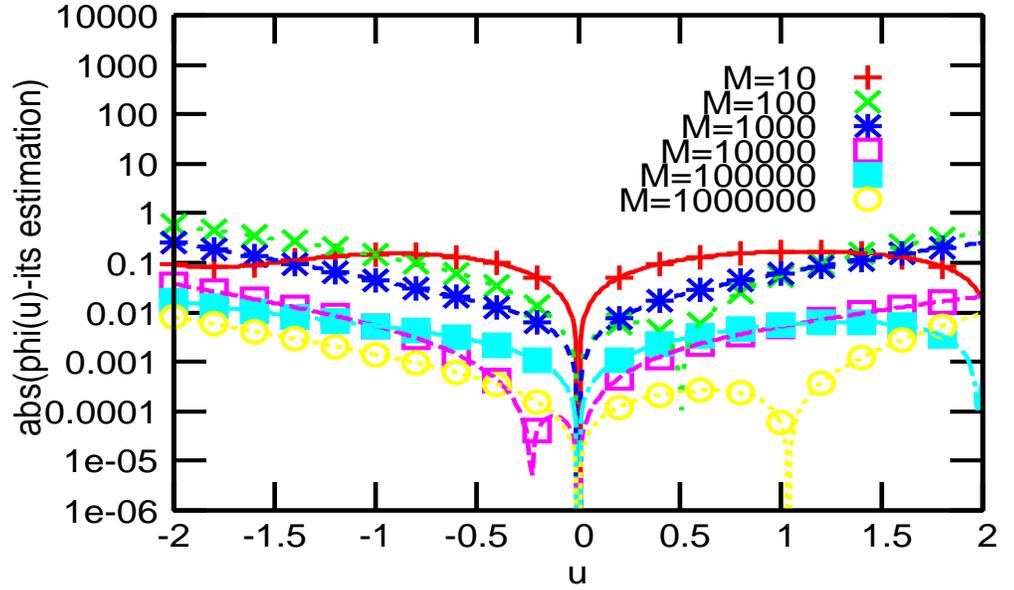} 
\caption{\label{fig1-est}
Relationship between the argument of the cumulant generating function 
when the demand $D_a$ is i.i.d.\ with 
a standard normalized distribution, with $L=1$ and
 $\left|\phi(u)-\log\left(\f{1}{M}\sum_{a=1}^Me^{uD_a}\right)\right|$.
The horizontal axis shows 
$u$, and 
the vertical axis shows 
the absolute difference between the 
cumulant generating function $\phi(u)$ and 
its estimate as  $\log\left(\f{1}{M}\sum_{a=1}^Me^{uD_a}\right)$, 
that is, 
 $\left|\phi(u)-\log\left(\f{1}{M}\sum_{a=1}^Me^{uD_a}\right)\right|$.
 When $M$ is large, 
$\log\left(\f{1}{M}\sum_{a=1}^Me^{uD_a}\right)$
is close to the  cumulant generating function $\phi(u)$.
Moreover, it can be verified that this result has the self-averaging property \cite{Shinzato-SA,Wakai}.
}
\end{center}
\end{figure}

\section{Analytic Results and Numerical Experiments}
In this section, 
we verify the effectiveness of our proposed method 
using several examples of a
 multidimensional normal distribution, 
which can be used to model correlated commodities for which the statistical properties are 
well known.
\subsection{Single commodity with i.i.d.\ normal distribution}
First, we define the difference of the demand at $t$, $D_t$, of a single commodity over lead time $L$ from the mean $\mu=E[D_t]$ as 
$X_t(=D_t-\mu)$; we assume that the novel variable $X_t$ is 
i.i.d.\ with the normal distribution 
with mean $E[X_t]=0$ and variance $E[X_t^2]=\s^2$. 
From this, we have $\phi(u)=\f{\s^2u^2}{2}$ 
when the safety stock quantity is $L\ve$, and we obtain 
\bea
\label{eq31-p}
Pr\left[L\ve\le X\right]\eq\int_{L\ve}^\infty dXp(X)\nn
&\le&e^{-L\left(\f{\ve^2}{2\s^2}\right),}
\eea
where $X=\sum_{t=1}^LX_t=\sum_{t=1}^LD_t-L\mu$. 
From this, 
when the allowable stockout rate is $\d$, the
safety stock quantity is $L\ve=\sqrt{-2L\s^2\log\d}$.

One point should be noted here. 
In the case of the normal distribution, 
since the cumulant generating function $\phi(u)=\f{\s^2u^2}{2}$ is already known, 
although the variance $\s^2$ appears in the safety stock quantity 
obtained from the above discussion, 
with respect to the demand distribution in general,
the variance does not always appear in the description of the cumulant 
generation function, 
and it is not possible to determine 
a safety stock quantity for the general case 
that can guarantee a more adequate level of stockout rate.
In a way similar to that shown in \hyou{fig1-est},
if we can estimate the cumulant generating function, 
we can use the rate function 
to more accurately determine the safety stock quantity; 
this practical stock management strategy 
uses our proposed method.

Moreover, note that for this model, $Pr[L\ve\le X]$ can be obtained directly, as follows:
\bea
\label{eq32-p}Pr\left[L\ve\le X\right]=H\left(\f{L\ve}{\sqrt{L\s^2}}
\right),
\eea
where
\bea
H(k)\eq\int_k^\infty 
\f{dt}{\sqrt{2\pi}}e^{-\f{t^2}{2}}.
\eea
From this and \siki{eq32-p}, when 
the allowable stockout rate is $\d$,
the safety stock quantity $L\ve (=SS_{\rm pre.})$ 
is 
\bea
SS_{\rm pre.}\eq\sqrt{L\s^2}H^{-1}(\d),
\eea
where 
$k=H^{-1}(\d)$ is the inverse function of $H(k)=\d$. It can be seen that 
the previous stock management strategy is a special case of our proposed method.

Two points should be noted here. First, 
although 
we analytically derive the safety stock quantity 
for the case of a normal distribution, 
analyze the Chernoff inequality in \siki{eq31-p}, and rigorously 
describe \siki{eq32-p}, 
it is not always easy to 
analytically solve for
the stockout rate $Pr[L\ve\le X]$. Next, 
in general, 
even if the probability density function, $P(X_t)$, of the demand difference from its mean, 
$X_t$,  is well known,
a rigorous allowable stockout rate $Pr[L\ve\le 
X]=E[\Theta(\sum_{t=1}^LX_t-L\ve)]$ is analytically described as follows: 
\bea
&&Pr[L\ve\le X]\nn
\eq
\area \prod_{t=1}^LdX_tP(X_t)\Theta\left(\sum_{t=1}^LX_t-L\ve\right)\nn
\eq\area \prod_{t=2}^LdX_tP(X_t)\nn
&&\qquad\times\int_{{L\ve-\sum_{t=2}^LX_t}}^\infty dX_1P(X_1).
\qquad\label{eq13}
\eea
However, 
even if we are able to find a closed form for the integral of $X_1$ 
in \siki{eq13},
the integrals of 
$X_2,\cdots,X_L$ are not easy to evaluate, and
we need to approximate this convolution integral, such as by using trapezoidal integration.
In particular, since it is necessary to compute the exponential of $L-1$, it cannot be rigorously evaluated. 
On the other hand, if we can examine the moment-generating function using the approach based on Chernoff inequality, then 
we can assess the rate function $R(\ve)$ and the safety stock quantity $L\ve$ corresponding to the allowable stockout rate.
\subsection{Definition of stockout rate}
In this subsection, 
we summarize our method for management of the safety stock of multiple commodities. 
For simplicity, we will consider the case of two types of commodities. For a lead time $L$, let $X_t$ and $Y_t$ be the normalized demands (that is, the difference between current demand and mean demand) of two commodities, respectively.
We assume that $X_t,Y_t$ are i.i.d.\ with the Gaussian distribution such that 
$E[X_t]=E[Y_t]=0$, $E[X_t^2]=\s_X^2$, $E[Y_t^2]=\s_Y^2$, and $E[X_tY_t]=\rho\s_X\s_Y$. From this, the probabilities that 
the sums of demand over lead time $L$, $X=\sum_{t=1}^LX_t$ and $Y=\sum_{t=1}^LY_t$, exceed $L\ve_X$ and $L\ve_Y$, 
respectively, are calculated as follows:
\bea
Pr[L\ve_X\le X]\eq H\left(\f{L\ve_X}{\sqrt{L\s^2_X}}\right),\\
Pr[L\ve_Y\le Y]\eq H\left(\f{L\ve_Y}{\sqrt{L\s^2_Y}}\right).
\eea
In addition, it is also necessary to evaluate the stockout probability: $Pr[L\ve_X\le X,L\ve_Y\le 
Y]$ and  $Pr[L\ve_X\le X,L\ve_Y\le 
Y]+Pr[L\ve_X>X,L\ve_Y\le Y]+Pr[L\ve_X\le X,L\ve_Y>Y]$
in the case of more than three commodities.
Since the variation of these states of stockout is exponentially increasing with the number of commodities, and 
from the discussion in subsection 
\ref{sec2.2}, the Chernoff inequality holds for any pattern of stockout. 
Hereafter, we will discuss just one of the stockout patterns, the case of all commodities are stockout.
\subsection{
Two commodities that are i.i.d.\ with normal distribution
\label{sec3.3}}
Let the demands of commodities 1 and 2 be
$D_t^X$ and $D_t^Y$, respectively.
Differences from the demand and its average demand over the lead time 
($E[D_t^X]=\mu_X,E[D_t^Y]=\mu_Y$) are represented as 
$X_t=D_t^X-\mu_X,Y_t=D_t^Y-\mu_Y$, 
$E[X_t]=E[Y_t]=0$, $E[X_t^2]=\s_X^2,E[Y_t^2]=\s_Y^2$, and 
$E[X_tY_t]=E[X_t]E[Y_t]=0$. 
From this, when the respective safety stock quantities are $L\ve_X$ and $L\ve_Y$ 
over lead time $L$, the stockout rate $Pr[L\ve_X\le X,L\ve_Y\le Y]$, which is the probability that both commodities are missing simultaneously, is 
\bea
&&
Pr\left[L\ve_X\le X,L\ve_Y\le Y\right]
\nn
\eq\int_{L\ve_X}^\infty dX\int_{L\ve_Y}^\infty dYp(X,Y)\nn
&\le&e^{-L\left(\f{\ve_X^2}{2\s^2_X}+\f{\ve_Y^2}{2\s^2_Y}\right)},
\eea
where 
$X=\sum_{t=1}^LX_t$ and $Y=\sum_{t=1}^LY_t$. From this, the allowable stockout rate when the 
safety stock quantities $L\ve_X,L\ve_Y$ are {
$L\ve_X=\sqrt{2L\s_X^2\log(1/\d_X)},L\ve_Y=\sqrt{2L\s_Y^2\log(1/\d_Y)}$} is 
$\d=\d_X\d_Y$.
Moreover, 
$Pr[L\ve_X\le X,L\ve_Y\le Y]$ can be directly evaluated as follows:
\bea
&&Pr\left[L\ve_X\le X,L\ve_Y\le Y\right]\nn
\eq 
H\left(\f{L\ve_X}{\sqrt{L\s^2_X}}\right)
H\left(\f{L\ve_Y}{\sqrt{L\s^2_Y}}\right)\label{eq16}.
\eea
From this, if we have
\bea
L\ve_X\eq\sqrt{L\s^2_X}H^{-1}\left(\d_X\right),\label{eq26}\\
L\ve_Y\eq\sqrt{L\s^2_Y}H^{-1}\left(\d_Y\right),\label{eq26-2}
\eea
then the allowable stockout rate $\d=\d_X\d_Y$ is obtained.

\subsection{Two normally distributed, correlated commodities \label{sec3.4}}
Let $X_t,Y_t$ be the demands of two correlated commodities over lead time $L$, where they are normally distributed with mean 
$E[X_t]=E[Y_t]=0$, variance $E[X_t^2]=\s_X^2$, $E[Y_t^2]=\s_Y^2$, and covariance $E[X_tY_t]={\rho}\s_X\s_Y$. Then, 
when the safety stock quantities are $L\ve_X$ and $L\ve_Y$, 
the stockout rate 
$Pr[L\ve_X\le X,L\ve_Y\le Y]$ is  
\bea
&&
\label{eq19}
Pr\left[L\ve_X\le X,L\ve_Y\le Y\right]\nn
\eq
\int_{L\ve_X}^\infty dX\int_{L\ve_Y}^\infty dYp(X,Y)\nn
&\le&\exp\left[-\f{L\left(\f{\ve_X^2}{\s^2_X}-2\rho\f{\ve_X\ve_Y}{\s_X\s_Y}+\f{\ve_Y^2}{\s_Y^2}\right)}{2(1-\rho^2)}\right],
\eea
where $X=\sum_{t=1}^LX_t$ and $Y=\sum_{t=1}^LY_t$.
{From this,
\bea
L\ve_X\eq\sqrt{-2L\s_X^2\log\d_X}\label{eq28},\\
L\ve_Y\eq\f{\s_Y}{\s_X}\rho(L\ve_X)
+\sqrt{-2L\s_Y^2(1-\rho^2)\log\d_Y
}\label{eq28-2},
\eea
and thus we can evaluate $\d=\d_X\d_Y$.}
Note that in this model, $Pr[L\ve_X\le X,L\ve_Y\le Y]$ is directly calculated 
as a double integral:
\bea
&&Pr\left[L\ve_X\le X,L\ve_Y\le Y\right]\nn
\eq \int_{\f{L\ve_X}{\sqrt{L\s_X^2}}
}^\infty 
\f{du}{\sqrt{2\pi}}e^{-\f{u^2}{2}}\nn
&&\times 
H\left(\f{1}{\sqrt{1-\rho^2}}\left(\f{L\ve_Y}{\sqrt{L\s^2_Y}}-\rho u\right)
\right),
\qquad
\label{eq29}
\eea
where, if $\rho=0$, \siki{eq29} matches \siki{eq16}.
In this case, we have another analytic description of the stockout rate: 
\bea
&&Pr[L\ve_X\le 
X,L\ve_Y\le Y]\nn
\eq\int_{\f{L\ve_X}{\sqrt{L\s_X^2(1-\rho^2)}}}^\infty 
ds
\int_{\f{L\ve_Y}{\sqrt{L\s_Y^2(1-\rho^2)}}}^\infty dt\f{1}{{2\pi}}e^{-\f{s^2+t^2}{2}}
\nn
&&
e^{\rho st}
\sqrt{1-\rho^2}.\qquad\label{eq30}
\eea
As in the discussion of \siki{eq13}, as the lead time $L$ increases, 
the calculation amount increases exponentially with the number of items to be considered ($N$), 
so for a sufficiently large $N$, it is not practical to use a direct method, such as Eq. (\ref{eq29}) or (\ref{eq30}).

\subsection{Numerical experiments}

\begin{figure}[bt]
\begin{center}
\includegraphics[width=1.0\hsize,height=0.6\hsize]{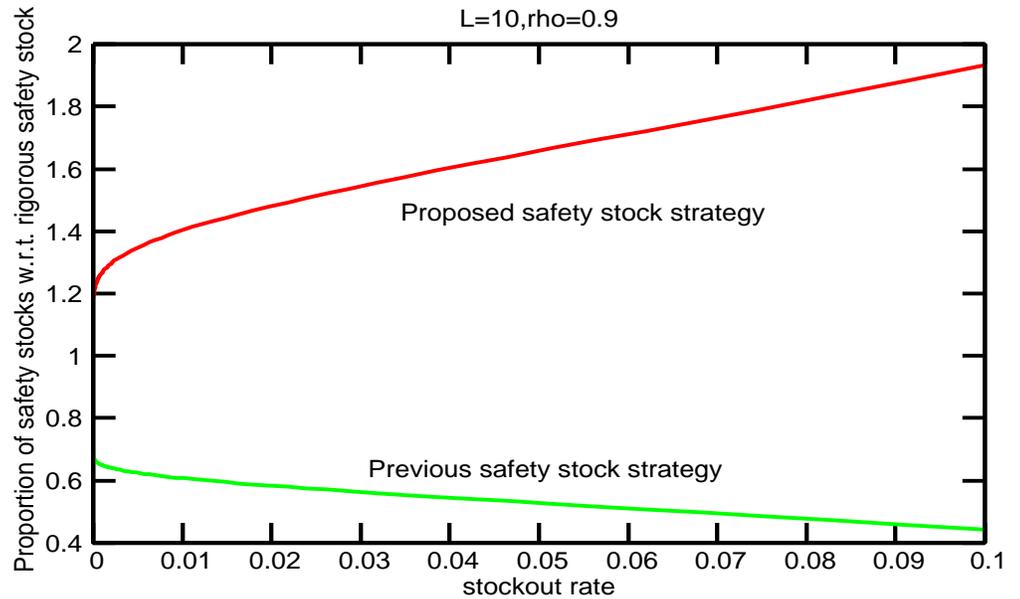} 
\caption{\label{fig1}
Comparison between the safety stock quantity derived by our proposed method 
 $SS_{\rm pro.}$ relative to the 
rigorous safety stock $SS_{\rm 
rig.}$  and  the safety stock quantity estimated by the previous method $SS_{\rm pre.}$ relative to the
rigorous safety stock $SS_{\rm 
rig.}$. The horizontal axis shows the allowable stockout rate $\d$ and 
the vertical axis shows the safety stock rate.
}
\end{center}
\end{figure}

\begin{figure}[bt]
\begin{center}
\includegraphics[width=1.0\hsize,height=0.6\hsize]{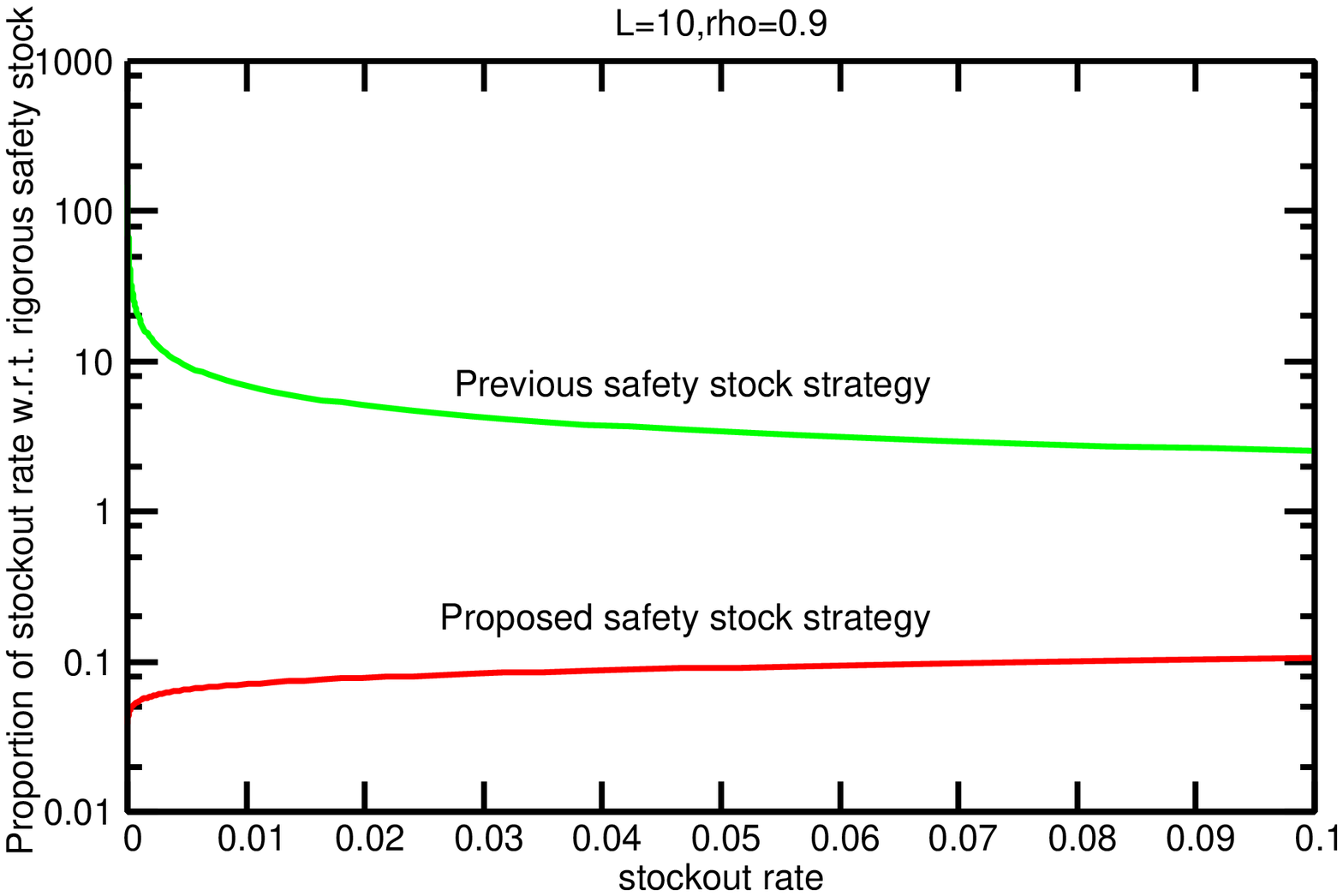} 
\caption{\label{fig2}
Comparison between
the proportion of the stockout rate by our proposed method with respect to the allowable stockout rate 
${\cal P}(SS_{\rm pro.})/{\cal P}(SS_{\rm rig.})$ and 
theproportion of the stockout rate by the previous method with respect to the allowable stockout rate 
${\cal P}(SS_{\rm pre.})/{\cal P}(SS_{\rm rig.})$.
The horizontal axis shows the allowable stockout rate $\d$ and the
vertical axis shows the proportionss of the stockout rate with respect to allowable stockout rate.
}
\end{center}
\end{figure}

In this subsection, we evaluate the effectiveness of the proposed method by conducting numerical experiments. 
For simplicity, we consider the case in which the demands have equal variance,  
$\s_X^2=\s_Y^2=\s^2$. From this symmetry, we can evaluate 
$\d_X,\d_Y$ by assuming  $ L\ve_X=L\ve_Y=L\ve$.
From Eqs. (\ref{eq26}), (\ref{eq26-2}), 
(\ref{eq28}), (\ref{eq28-2}) and (\ref{eq30}),
for a given allowable stockout rate
$\d$, 
the rigorous safety stock quantity $SS_{\rm rig.}$ (rig.=rigorous),
the proposed safety stock quantity $SS_{\rm pro.}$ (pro.=proposed), and 
the previous safety stock quantity $SS_{\rm pre.}$ (pre.=previous) can be evaluated:
\bea
SS_{\rm rig.}\eq{\cal P}^{-1}(\d),\\
SS_{\rm pro.}\eq\sqrt{{L}
\s^2({1+\rho})\log(1/\d)},\\
SS_{\rm pre.}\eq\sqrt{L\s^2}H^{-1}(\sqrt{\d}),
\eea
where we use the inverse function of ${\cal P}(SS)=Pr[SS\le X,SS\le Y]$ defined in Eqs. (\ref{eq29}) and (\ref{eq30}).
Furthermore, from subsection \ref{sec3.3}, we have
$\d_X=\d_Y=\sqrt{\d}$, and from subsection  
\ref{sec3.4}, we have $\d_X=\d^{\f{1+\rho}{2}},\d_Y=\d^{\f{1-\rho}{2}}$.

\hyou{fig1} compares
the proportion of the proposed safety stock quantity with respect to the rigorous safety stock quantity 
$SS_{\rm pro.}/SS_{\rm 
rig.}$ and 
the proportion of the previous safety stock quantity with respect to the rigorous safety stock quantity 
$SS_{\rm pre.}/SS_{\rm rig.}$. Similarly, 
the proportions of the stockout rates 
${\cal P}(SS_{\rm pro.})/{\cal P}(SS_{\rm rig.})
$ and ${\cal P}(SS_{\rm pre.})/{\cal P}(SS_{\rm rig.})$ are indicated in 
\hyou{fig2}.
In Figs. \ref{fig1} and \ref{fig2}, $\s^2=1$, $\rho=0.9$,
$L=10$.
From \hyou{fig1}, it can be seen that 
when the allowable stockout rate is less than $10\%$, 
the safety stock quantity derived from the proposed method 
is 1.2 to 1.9 times the rigorous safety stock quantity.
From \hyou{fig2}, it can be seen that
the stockout rate of the proposed method is always less than the allowable stockout rate; this is guaranteed by the Chernoff inequality.
However, with the existing method, the stockout rate always exceeds the allowable value.
In addition, 
 the other cases of $\rho,L$ used in Figs. \ref{fig1} and \ref{fig2} are also guaranteed. 
Thus, it is clear that despite the safety stock, which is meant to minimize 
the opportunity loss, 
when there is a correlation between the different types of commodities,
the proper level cannot be guaranteed by the previous method.

\section{Conclusion and future work}

In this study, we proposed a novel method for stock management of multiple commodities with correlated demands; we used the Chernoff inequality to derive the relation 
between the safety stock quantity and the allowable shortage rate. 
Using the rate function defined by Legendre's transformation of the cumulant generating function, we were able to determine the safety stock level that would ensure  
that the stockout rate would always be acceptable.
The theoretical properties are well known, since they follow a multidimensional normal distribution. We derived the safety stock for three different models: our proposed method, an existing method, and the exact (analytical) method. We determined the allowable stockout rate for each, and we verified the effectiveness of our method by conducting numerical experiments.

In this study, we considered two correlated commodities, but in order to increase the usefulness of the proposed method, 
we discussed the safety stock amount for multiple commodities, including those with positive and negative correlations. 
We discussed the safety stock amount for one lifecycle, but in reality, it is important to consider the safety stock while including   
preorders and back orders and spanning multiple life cycles. 
We intend to address this in our future research.

\section*{Acknowledgements}
The author appreciates the fruitful comments of I.\ Arizono, Y.\ Takemoto, K.\
Kobayashi, and I.\ Kaku. 
This work was supported in part
by Grant-in-Aid No.\ 15K20999; 
the President Project for Young Scientists at Akita Prefectural University; 
Research Project No.\ 50 of the National Institute of Informatics, Japan; 
Research Project No.\ 5 of the Japan Institute of Life Insurance; 
Research Project of the Institute of Economic Research Foundation at Kyoto University;
Research Project No.\ 1414 of the Zengin Foundation for Studies in Economics and Finance; 
Research Project No.\ 2068 of the Institute of Statistical Mathematics; 
Research Project No.\ 2 of the Kampo Foundation; 
and Research Project of the Mitsubishi UFJ Trust Scholarship Foundation.

\appendix

\section{Properties of the rate function\label{app-a1}}
In this appendix, we present the properties of the rate function $R(\eta)$, which is used to derive the relation 
between the safety stock quantity and the allowable
stockout rate. First, we can rewrite \siki{eq10}, the definition of the rate
function, as follows:
\bea
\label{eq20}
R(u,\eta)
\eq
u\eta-
\log E[e^{uD}],
\eea
where $R(0,\eta)=0$, and $R(u,\eta)$ is a continuous function of $u$, and so
\bea
\label{eq21}
R(\eta)&\ge&
0.
\eea
Thus, we obtain a positive upper bound for
$Pr[\eta\le D]$ 
by using the Chernoff inequality: 
 $Pr[\eta\le D]\le e^{-R(\eta)}\le e^{-0}=1$.

Next, since the exponential function is convex,
$E[e^{uD}]\ge e^{uE[D]}$ holds, and so we obtain
\bea
\label{eq22}
R(u,\eta)
&\le&
u\eta- uE[D].
\eea
When $\eta\le E[D]$, since $u(\eta-E[D])\le0$, the
optimal value for $u$ is close to $0$, and we obtain $R(\eta)=0$. On the other
hand, when $\eta\ge E[D]$, 
the optimal value for $u$ is not equal to $0$, and so we obtain
$R(\eta)>0$. Thus, we can rewrite the definition for $\eta$ as $\eta=E[D]+\ve$, 
and it can be seen that the sign of the rate function is determined by the sign of $\ve$,
that is, if $\ve\le0$, then $R(\eta)=0$, otherwise, $R(\ve)>0$. 
Hereafter,
we will consider the rate function $R(\eta)$ to be a function of the distance from the
mean, $E[D]$.

We also note that, from the definition
of the rate function, $R(u,\eta)$ is convex with respect to $u$, since
 $\pp{^2R(u,\eta)}{u^2}\le0$.
 Thus, the following updating rule will optimize $u$:
\bea
\label{eq23}
u_{s+1}
\eq
u_s+
\kappa_u
\lim_{u\to u_s}
\pp{R(u,\eta)}{u},
\eea
where $u_s$ is the state of $u$ at step $s$, for $s\in{\bf Z}$, and the step constant is $\kappa_u=10^{-3}$.  
This iterative definition ceases if $|u_s-u_{s+1}|$ is close to $0$ 
or
the variation of $u$ is smaller than an infinitesimal $\theta$,
 such as $\theta=10^{-6}$.
Both the control
parameter $u$ 
and the constant $\eta$ can develop to represent a multidimensional case;
for instance, for $\vec{u},\vec{\eta}\in{\bf R}^N$, 
$R(\vec{u},\vec{\eta})$ 
is convex 
with respect to $\vec{u}$, and thus 
$R(\vec{\eta})=\mathop{\max}_{\vec{u}}
R(\vec{u},\vec{\eta})$
is concave with respect to $\vec{\eta}$.
From the above conclusions, and using the convexity of the rate function, the duality of the
rate function and the cumulant generating function, and an algorithm developed in convex optimization research,
we can assess a more 
appropriate safety stock quantity 
for stock management under stochastic phenomena.

\section{Availability of \siki{eq10}\label{appendix:A}}
\siki{eq10} shows the Chernoff inequality with respect to $N$ commodities. 
The normalized demand (or the difference between current demand and mean demand) over lead time $L$, $X_{it},(i=1,\cdots,N,t=1,\cdots,L)$, 
is not temporally correlated and follows the normal distribution with 
$E[X_{it}]=0,E[X_{it}X_{jt}]=\Sigma_{ij}$. From \siki{eq10}, we can derive the following inequality:
\bea
Pr\left[L\vec{\ve}\le\vec{X}\right]&\le&e^{-L\vec{u}^{\rm 
T}\vec{\ve}+\f{L}{2}\vec{u}^{\rm T}\Sigma\vec{u}},
\eea
where $\vec{u}\ge0$, 
$\vec{X}=\sum_{t=1}^L\vec{X}_t$, and the variance-covariance matrix $\Sigma=\left\{\Sigma_{ij}\right\}\in{\bf R}^{N\times N}$. We can obtain a tighter upper bound by using the following equation with respect to the parameter 
$\vec{u}$:
\bea
\vec{u}_{s+1}\eq\max\left(\vec{u}_s+\kappa(\vec{\ve}-\Sigma\vec{u}_s),0\right),
\eea
where $\kappa$ is positive and infinitesimal, $\vec{u}_0=(1,1,\cdots,1)^{\rm T}\in{\bf R}^N$, and at step $s$, 
$\vec{u}_s=(u_{1,s},\cdots,u_{N,s})^{\rm T}\in{\bf R}^N$. The
stopping condition is $\Delta=\sum_{i=1}^L\left|u_{i,s+1}-u_{i,s}\right|<10^{-6}$. The rigorous stockout rate $L\vec{\ve}$ is calculated as follows:
\bea
Pr\left[L\vec{\ve}\le\vec{X}\right]\eq
\int_{L\vec{\ve}}^\infty \f{d\vec{X}e^{-\f{1}{2L}\vec{X}^{\rm T}\Sigma^{-1}\vec{X}}}{(2\pi)^{\f{N}{2}}\sqrt{\det(L\Sigma)}}.
\qquad\label{eq20}
\eea
We first diagonalize the variance covariance matrix $\Sigma$, and then we prepare 
the novel variables, which are in the directions of the eigenvectors of $\Sigma$.
However, 
since the integral domain of the novel variables is more complicated, 
it is not easy to analytically assess the right-hand side of \siki{eq20}.
Thus, we use trapezoidal integration to approximate it when 
the number of commodities, $N$, is large; however,
the computational complexity increases 
exponentially with $N$, so
this approach is not practical.

\section{Safety stock management for fungible commodities \label{app2}}
We here discuss the 
safety stock management for 
the case in which $X$ and 
$Y$ are fungible commodities.
Using the model setting in subsection \ref{sec3.4}, 
with respect to the modified demands for two commodities $X$ and $Y$ over lead time $L$, where we consider 
their differences from their means, 
$X_t,Y_t$ (since they are fungible, we assume that they have a negative correlation), 
the probability that 
the sum 
$\sum_{t=1}^LX_t+\sum_{t=1}^LY_t$ is larger than 
the safety stock quantity $L\ve$ is calculated using the Chernoff inequality:
\bea
&&Pr\left[L\ve\le\sum_{t=1}^L(X_t+Y_t)\right]\nn
&\le& 
e^{-uL\ve}E\left[e^{u\sum_{t=1}^L(X_t+Y_t)}\right]\nn
\eq
e^{-uL\ve+\f{Lu^2}{2}(\s_X^2+\s_Y^2+2\rho\s_X\s_Y)}.
\eea
From this, the rate function is
\bea
R(\ve)\eq\mathop{\max}_{u>0}
\left\{
u\ve-\f{u^2}{2}(\s_X^2+\s_Y^2+2\rho\s_X\s_Y)
\right\}\nn
\eq\f{\ve^2}{2(\s_X^2+\s_Y^2+2\rho\s_X\s_Y)}.
\eea
Thus, the safety stock with respect to the allowable stockout rate $\d$ is
\bea
L\ve\eq\sqrt{-2L(\s_X^2+\s_Y^2+2\rho\s_X\s_Y)\log\d}.\label{eq32}
\eea
Here, we assume that the correlation between $X$ and $Y$ is negative, but 
this also 
applies when the correlation is positive. Thus, 
the safety stock quantity $L\ve$
is a monotonically nondecreasing function with respect to the 
correlation coefficient $\rho$.


\bibliographystyle{amsplain}

\end{document}